\documentstyle[12pt]{article}
\normalbaselines
\begin{document}
\newcommand{\C}{{\bf C}}
\newcommand{\Z}{{\bf Z}}
\newcommand{\R}{{\bf R}}

\title{Congruences for real algebraic curves on an ellipsoid}
\author{Grigory Mikhalkin}
\date{}
\maketitle

\begin{abstract}
The problem of arrangement of a real algebraic curve on a real algebraic
surface is related to the 16th Hilbert problem.
We prove in this paper new restrictions on arrangement of nonsingular real
algebraic curves on an ellipsoid.
These restrictions are analogues of Gudkov-Rokhlin, Gudkov-Krakhnov-Kharlamov,
Kharlamov-Marin congruences for plane curves (see e.g. \cite{V} or \cite{V1}).
To prove our results we follow Marin approach \cite{Marin} that is a
study of the quotient space of a surface under the complex conjugation.
Note that the Rokhlin approach \cite{R} that is a study of the 2-sheeted
covering of the surface branched along the curve can not be directly applied
for a proof of Theorem \ref{Rokhlin} since the homology class of a curve of
Theorem \ref{Rokhlin} can not be divided by 2 hence such a covering space
does not exist.
\end{abstract}

\section{Formulations of results}
It is well-known that a complex quadric is isomorphic to $\C P^1\times
\C P^1$ and an algebraic curve on a quadric is defined by a bihomogeneous
polynomial of bidegree $(d,r)$.
If the curve is real and the quadric is an ellipsoid then $d=r$
and the curve can be represented as the intersection of an ellipsoid
and a surface of degree $d$ in $\R P^3$ (this is because a curve
of bidegree $(d,r), d\neq r$ can not be invariant under the involution
of the complex conjugation of ellipsoid).

Let $Q$ be an ellipsoid, $\R Q$ and $\C Q$ be the spaces of its real and
complex points; $A$ be a nonsingular real algebraic curve of bidegree
$(d,d)$ on $Q$; $\R A$ and $\C A$ be the spaces of real and complex points
of $A$.
Components of $\R A$ are called ovals and the number of ovals of $\R A$
is denoted by $l$.
$\R A$ divides $\R Q$ into two parts with a common boundary.
Let $B_0$ and $B_1$ denote these parts in such a way that in the case
when $l$ is even congruence $\chi(B_0)\equiv 0\pmod{4}$ is correct
(as it is usual, $\chi$ denotes Euler characteristic).

According to F.Klein, $A$ is a curve of type I (II) if $\C A \setminus \R A$
is not connected (otherwise).
It is well-known (see \cite{R1}) that if $A$ is a curve of type I then
$\R A$ has two natural opposite orientations.

We need the following definitions to formulate Theorem \ref{Fiedler}.
Let us choose one of two complex orientation of $\R A$ and
some orientation of $B_0$.
Oval $C$ of $\R A$ is called disorienting if these orientations induce
opposite orientations on it.
$C$ divides $\R Q$ into two disks $D$ and $D'$.
Let $x(D)$ be equal to $\chi(B_1\cap D)\bmod{2}$.
It is clear that if $l\equiv 0\pmod{2}$ then $x(D)=x(D')$;
in this case we set $x(C)$ to be equal to $x(D)$.

The main results of this paper are the following.

\newtheorem{th}{Theorem}
\begin{th}
\label{Rokhlin}
Let $d$ be an odd number.
\begin{itemize}
\begin{description}
\item[a)] If $A$ is an M-curve (i.e. $l=(d-1)^2+1$) then
\begin{displaymath}
\chi(B_0)\equiv\chi(B_1)\equiv\frac{d^2+1}{2}\pmod{8}
\end{displaymath}
\item[b)] If $A$ is an (M-1)-curve (i.e. $l=(d-1)^2$) then
\begin{displaymath}
\chi(B_0)\equiv\frac{d^2-1}{2}\pmod{8}
\end{displaymath}
\begin{displaymath}
\chi(B_1)\equiv\frac{d^2+3}{2}\pmod{8}
\end{displaymath}
\item[c)] If $A$ is an (M-2)-curve (i.e. $l=(d-1)^2-1$) and
\begin{displaymath}
\chi(B_0)\equiv\frac{d^2-7}{2}\pmod{8}
\end{displaymath}
then $A$ is of type I.
\item[d)] If $A$ is of type I then
\begin{displaymath}
\chi(B_0)\equiv\chi(B_1)\equiv1\pmod{4}
\end{displaymath}
\end{description}
\end{itemize}
\end{th}

\begin{th}
\label{Fiedler}
Let $d$ be an even number.
\begin{itemize}
\begin{description}
\item[a)] \label{af} If $A$ is an M-curve and if all components of $B_1$ have
even Euler
characteristics then
\begin{displaymath}
\chi(B_0)\equiv d^2\pmod{16}
\end{displaymath}
\begin{displaymath}
\chi(B_1)\equiv 2-d^2\pmod{16}
\end{displaymath}
\item[b)] \label{bf} Let $A$ be a curve of type I with chosen complex
orientation.
If there exist an orientation of $B_0$ such that $x(C)=0$ for every
disorienting oval $C$ then
\begin{displaymath}
\chi(B_0)\equiv d^2\pmod{8}
\end{displaymath}
\begin{displaymath}
\chi(B_1)\equiv 2-d^2\pmod{8}
\end{displaymath}
\end{description}
\end{itemize}
\end{th}

{\bf\noindent Remark:} Propositions d) of Theorem \ref{Rokhlin}
and b) of Theorem \ref{Fiedler} can be deduced from the formula of complex
orientations of Zvonilov \cite{Z}.

\section{Proof of Theorem \ref{Rokhlin}}
Let $W$ denote $B_0\cup\C A/conj$.
As it was shown in \cite{L} a manifold $\C Q/conj$ where $conj$
is the involution of complex conjugation is diffeomorphic to
${\overline{\C P^2}}$
It is easy to check that
\begin{equation}
\label{WW}
W\circ W=d^2-2\chi(B_0)\equiv 1\pmod{2}
\end{equation}
Hence, $W$ is a characteristic surface in $\C Q/conj$.
Let us apply Guillou-Marin congruence \cite{GM}
\begin{displaymath}
\sigma(\C Q/conj)\equiv W\circ W+2\beta(q)\pmod{16}
\end{displaymath}
where $\sigma$ is signature, $\beta(q)$ is Brown invariant of quadratic
Guillou-Marin form \newline $q: H_1(W;\Z_2)\rightarrow \Z_4$ of surface $W$
in $\C Q/conj$.
This congruence, equality \ref{WW} and equality
$\sigma({\overline{\C P^2}})=-1$ follow that
\begin{displaymath}
\chi(B_0)\equiv\frac{d^2+1}{2}+\beta(q)\pmod{8}.
\end{displaymath}

The calculations of $\beta(q)$ below are similar to the calculations
in \cite{Marin}, \cite{KV}.
Let $L$ denote the subspace of $H_1(W;\Z_2)$ generated by classes realized
by ovals of $\R A$.
Let $U=L^{\perp}$ denote the orthogonal complement of $L$
with respect to the intersection form of $W$.
It is clear that $U^{\perp}=L\subset U$ and $q|_L=0$.
Therefore, $U$ is an informative subspace and $\beta(q)=\beta(q')$,
where $q'$ is the form on $U/U^{\perp}$ induced by $q$ (see \cite{KV}).
According to \cite{KV} one can get the following.
\begin{itemize}
\begin{description}
\item[a)] $dim(U/U^{\perp})=0$, hence $\beta(q')\equiv 0\pmod{8}$
\item[b)] $dim(U/U^{\perp})=1$, hence $\beta(q')\equiv\pm 1\pmod{8}$
\item[c)] $dim(U/U^{\perp})=2$, $\beta(q')\equiv 4\pmod{8}$, hence $q'$ is
even,
therefore $A$ is of type I
\item[d)] $q'$ is even, hence $\beta(q')\equiv 0\pmod{4}$.
\end{description}
\end{itemize}

\section{Lemma for Theorem 2}
Let $A_+$ be the closure of one of two components of $\C A\setminus\R A$.
Let $W_0$ denote $A_+\cup B_0$ and let $W_1$ denote $A_+\cup B_1$.
\newtheorem{lemma}{Lemma}
\begin{lemma}\label{ne0} Surface $W_1$ is not homologous to zero modulo 2 in
$\C Q$
\end{lemma}
{\em{\underline{Proof}}} If $W_1$ is homologous to zero then $W_1$ is a
characteristic surface in $\C Q$ and we can apply Guillou-Marin congruence
\begin{displaymath}
0\equiv W_1\circ W_1+2\beta(q_1)\pmod{16}
\end{displaymath}
where $q_1$ is the Guillou-Marin form of surface $W_1$ in $\C Q$.
But $W_1\circ W_1=d^2-\chi(B_1)$ and $\beta(q_1)\equiv 0\pmod{2}$
since $dim H_1(W_1;\Z_2)\equiv \chi(W_1)\equiv 0\pmod{2}$,
note that $\chi(W_1)$ is even, because $W_1$ is the result of gluing
of two orientable surfaces along whole boundary.
Congruences above and $d^2\equiv 0\pmod{4}$ imply that $\chi(B_1)\equiv 0
\pmod{4}$ follow that $\chi(B_1)\equiv 0\pmod{4}$.
This is a contradiction to the choice of $B_1$.

\section{Proof of Theorem 2}
Let $e_1$ and $e_2$ denote the elements of $H_2(\C Q)$ realized by generating
lines of quadric $\C Q$.
It is clear that $e_1$ and $e_2$ form a basis of $H_2(\C Q)$,
$conj_*e_1=-e_2$, $conj_*e_2=-e_1$.
Let $(\alpha,\beta)$ denote $\alpha e_1+\beta e_2, \alpha,\beta \in \Z$.
Each of generating lines transversally intersects $\R Q$ at one point,
hence $[\R Q]\equiv (1,1)\pmod{2}$.
{}From relations $[W_j]-conj_*[W_j]\equiv [\C A]\equiv (d,d)\pmod{2},j=0,1$,
$[W_0]+[W_1]\equiv [\R Q]\equiv (1,1)\pmod{2}$ and Lemma \ref{ne0}
one can deduce that $[W_0]\equiv 0\pmod{2}$.
Therefore $W_0$ is a characteristic surface in $\C Q$ and there is
Guillou-Marin form $q_0$ on $H_1(W_0;\Z_2)$.
Value of $q_0$ on the element of $H_1(W_0;\Z)$ realized by oval $C$
is equal to $x(C)$ (for the calculations see \cite{F}).
Let $L_0$ denote the subspace of $H_1(W_0;\Z_2)$ generated by such ovals C
that $x(C)=0$.
In the case \ref{af} we have $L_0=L_0^{\perp}$ therefore $\beta(q_0)\equiv 0
\pmod{8}$.
In the case \ref{bf} we have that form $q_0':L_0^{\perp}\rightarrow \Z_4$
induced by $q_0$ is even therefore $\beta(q_0)\equiv 0\pmod{4}$.
Now theorem 2 follows from Guillou-Marin congruence \cite{GM}
\begin{displaymath}
0\equiv d^2-\chi(B_0)+2\beta(q_0)\pmod{16}
\end{displaymath}
\section{Applications for the curves of bidegrees (3,3) and (5,5)}
Theorem \ref{Rokhlin} and Harnack inequality give a complete system
of restrictions for real schemes of flexible curves of bidegree (3,3)
(a definition of flexible curves on an ellipsoid is similar to the definition
given in \cite{V} for plane curves).
All real schemes realizable by flexible curves of bidegree (3,3) can be
realized by algebraic curves of bidegree (3,3)
(the classification of real schemes of algebraic curves of bidegree (3,3)
is $<$$\alpha$$>$, $\alpha\leq 5$ and $1\sqcup 1$$<$$1$$>$ -- see e.g.
\cite{Z}).

Real schemes of M-curves of bidegree (5,5) allowed by Theorem \ref{Rokhlin}
are $\alpha\sqcup 1$$<$$\beta$$>$$, \alpha+\beta=16,\beta\equiv 2\pmod{4}$ and
$\alpha\sqcup 1$$<$$\beta$$>$$\sqcup 1$$<$$\gamma$$>$$,\alpha+\beta+\gamma=15,
\alpha\equiv 1\pmod{4}$.
These schemes except might be $1\sqcup 1$$<$$6$$>$$\sqcup 1$$<$$8$$>$ and
$1\sqcup 1$$<$$5$$>$$\sqcup 1$$<$$9$$>$ can be constructed by smoothing
singularities of
images under birational transformations of appropriate plane M-curves
of degree 5 intersecting a real line at 5 different real points and
M-curves of degree 6 transversally intersecting a real line at 4 different real
points.

\section{An absence of congruences similar to Theorem 1
for curves of even bidegrees}
There is no congruence modulo 4 for each halves $B_1, B_2$ of a complement
of M-curve of even bidegree (because in this case $l\equiv 0\pmod{2}$,
hence $\chi(B_0)\equiv 0\pmod{2}$ and $\chi(B_0)\equiv\chi(B_1)+2\pmod{4}$).
A classification of the real schemes of curves of bidegree (4,4) shows
an absence of nontrivial congruences similar to Theorem 1 :
all schemes allowed by Harnack inequality and Bezout theorem are
realizable by algebraic curves - $<$$\alpha$$>$,$\alpha\leq 10$ and
$\alpha\sqcup 1$$<$$\beta$$>$, $\alpha+\beta\leq 9$
(see e.g.\cite{Z}).

Author is indebted to O.Ya.Viro, V.I.Zvonilov and V.M.Kharlamov for their
attention to the paper.

\appendix
\section{Congruences for curves with imaginary singularities}
This appendix can be regarded as the addendum to the paper of V.M.Kharlamov
and O.Ya.Viro \cite{KV}.
{}From theorems of paper \cite{KV} under appropriate birational transformations
one can deduce the congruences of Matsuoka \cite{Ma} for curves on a
hyperboloid.
But the congruences of these paper can not be deduced in this way,
because images of nonsingular curves on an ellipsoid under birational
isomorphism have imaginary singular points.
There are no consideration of curves with imaginary singular points in
\cite{KV}
, but the approach of \cite{KV} can be applied to these curves either.
Let us consider this application.

Link $L$ in $S^3$ is said to be even if the linking number of each component of
$L$ and the set of all other component is even.
The Arf-invariant of even link $L$ is equal by the definition to the
Arf-invariant of Seifert form of some Seifert surface of $L$.

Let $A$ be the plane real curve of degree $m=2k$ and let $\R P^2_+$ be one
of two parts of $\R P^2$ bounded by $\R A$.
We shall say that a singular point of $\C A$ is even if the link of this
singular point is even.
We shall set the Arf-invariant of singular point to be equal to the
Arf-invariant of the link of this singular point.

\newtheorem{prop}{}[section]
\begin{prop}
\label{a1}
Suppose that all singular points of $A$ are imaginary and even.
Let $Ar$ be the sum of Arf-invariants of singular points taken per one from
each
pair of complex-conjugated singular points.
Let $\R P^2_+$ be orientable.
\begin{itemize}
\begin{description}
\item[a)] If $A$ is an M-curve then $\chi(\R P^2_+)\equiv k^2+4Ar\pmod{8}$
\item[b)] If $A$ is an (M-1)-curve then $\chi(\R P^2_+)\equiv k^2\pm
1+4Ar\pmod{8}$
\item[c)] If $A$ is an (M-2)-curve of type II then $\chi(\R P^2_+)\equiv
k^2+d+4Ar\pmod{8}$, where $d\in\{0,2,-2\}$
\item[d)] If $A$ is of type I then $\chi(\R P^2_+)\equiv k^2\pmod{4}$
\end{description}
\end{itemize}
\end{prop}

The proof of Theorem \ref{a1} is similar to the proofs of theorems in
\cite{KV}.

Theorem \ref{Rokhlin} of present paper can be deduced from the application
of \ref{a1} to plane curves of degree $2d$ with two imaginary
complex-conjugated nondegenerate singular points of multiplicity $d$ with the
help of \ref{a4}
(it is easy to see that a nondegenerate singular point of multiplicity $d$ is
even iff $d$ is odd).

Let us extend theorem 3.A of \cite{KV} to the case of singular curves with
imaginary singular points (for the notations in \ref{a2} see \cite{KV}).

\begin{prop}
\label{a2}
Suppose that $\Z_4$-quadratic form $\tilde{Q_{\Delta}}$ is informative
and suppose that all imaginary singular points of $A$ are even.
\begin{itemize}
\begin{description}
\item[a)] If $A$ is an M-curve then $\chi(\R P^2_+)\equiv
k^2+\beta(\tilde{q_{\Delta}})+\tilde{b}+4Ar\pmod{8}$
\item[b)] If $A$ is an (M-1)-curve then $\chi(\R P^2_+)\equiv k^2\pm
1+\beta(\tilde{q_{\Delta}})+\tilde{b}+4Ar\pmod{8}$
\item[c)] If $A$ is an (M-2)-curve of type II then $\chi(\R P^2_+)\equiv
k^2+d+\beta(\tilde{q_{\Delta}})+\tilde{b}+4Ar\pmod{8}$, where $d\in\{0,2,-2\}$
\item[d)] If $A$ is of type I then $\chi(\R P^2_+)\equiv
k^2+\beta(\tilde{q_{\Delta}})+\tilde{b}\pmod{4}$
\end{description}
\end{itemize}
\end{prop}

\newtheorem{ex}{Example}
\begin{ex}
Let us consider the application of \ref{a2} to the problem of apparent contour
of real cubic surface.
Let $\R B$ be a surface of degree 3 in $\R P^3$.
Point $a\in \R P^3$ determine projection $pr:\R B\rightarrow\R P^2$.
The image of singular point of this projection is called the apparent contour
of cubic.
A generic apparent contour of cubic is a real plane singular curve with 6
singular points -- cusp points (this can be deduced from the fact that
an apparent contour is determined by the discriminant of cubic polynomial).
\end{ex}

\begin{prop}
\label{a3}
An apparent contour of cubic surface can not have real schemes represented on
figures 1 and 2
\end{prop}

{\underline{\em{Proof}}} It is easy to check that the Arf-invariant of cusp
point is equal to 1.
Curves on figures 1 and 2 have one pair of imaginary
complex-conjugated cusp points.
It is easy to calculate $\beta(\tilde{q_{\Delta}})$ with the help of
perturbation of curve.
After the calculation we have a contradiction with \ref{a2}.

We shall say that multiplicitive sequence of singular point of curve
(for the definition see \cite{Br}) is odd if the sum of numbers in every
round in sequence is odd (i.e. the intersection number of
the exceptional divisor of every
$\sigma$-process of resolution of singularity and the proper preimage of
the curve under this $\sigma$-process is odd).
Let $s_j$ denote this sum, where $j$ is the number of $\sigma$-process.
It is easy to see that singular points with odd multiplicitive sequence
are even.
Let us calculate Arf-invariants of these points.

\begin{prop}
\label{a4}
The Arf-invariant of singular points with odd multiplicitive sequence
of curve $\C A$ is equal to $\sum_j\frac{s^2_j-1}{8}$
\end{prop}

{\underline{\em{Proof}}} Let us glue with the evident diffeomorphism
of boundaries the pair $(D^4,D^4\cap\C A)$ and the pair $(D^4,D^4\cap\C
A_{\epsilon})$, where $D^4$ is a small ball with the center in the singular
point
and $A_{\epsilon}$ is a very small perturbation of $A$.
Let $(X,C)$ denote the result of gluing.
$C$ realizes the homology class dual to $w_2(X)$, because $X\approx S^4$.
Let us resolve by $\sigma$-processes the singular point of $(X,C)$.
Let $(M,S)$ denote the result of resolution.
It is easy to see that the self-intersection number of surface $S$ in
manifold $M$ is equal to $\sum_j (-s^2_j)$ and the signature of $M$
is equal to $\sum_j (-1)$ (i.e. to the number of $\sigma$-processes
multiplied by (-1) ).
The fact that $S$ is a characteristic surface in $M$ follows from the fact that
the multiplicitive sequence of the singular point is odd.
It is easy to see that the Arf-invariant of surface $S$ is equal to
the Arf-invariant of the singular point.
Now \ref{a4} follows from Rokhlin congruence
\begin{displaymath}
8Arf(M,S)\equiv\sum_j (1-s^2_j)\pmod{16}
\end{displaymath}

\section{Further generalizations}

In this appendix we announce results generalizing Theorem \ref{Rokhlin},
Rokhlin congruence \cite{R} and Kharlamov\cite{Kh}-Gudkov-Krakhnov\cite{GK}
congruence for curves on surfaces.
Proofs of these results will be published separately.

First consider the absolute case.

Let $\R B$ be a nonsingular real algebraic surface,
$\C B$ be its complexification.
Let $D:H^*(\C B;\Z_2)\rightarrow H_*(\C B;\Z_2)$ be the operator of
Poincar\'{e} duality.

\begin{prop}
\label{b1}
If $Dw_2(\C B)=[\R B]$ then there exists a natural separation of $\R B$
into two closed surfaces $B_1$ and $B_2$.
This separation is determined by the condition that $B_j, j=1,2$ is
a characteristic surface in $\C B/conj$.
There is congruence
\begin{displaymath}
\chi(B_j)\equiv\frac{\chi(\R B)-\sigma(\C B)}{4}+
\beta(q|_{H_1(B_j;\Z_2)})\pmod{8}
\end{displaymath}
where q is the Guillou-Marin form of surface $\R B$ in $\C B$
\end{prop}

We shall call surfaces $B_1$ and $B_2$ the surfaces of complex separation.

\begin{prop}[Corollary]
\label{b2}
If $Dw_2(\C B)=[\R B]$ and $\R B$ is connected surface then
\begin{displaymath}
\chi(\R B)\equiv\sigma(\C B)\pmod{32}
\end{displaymath}
\end{prop}

\begin{prop}[Corollary]
\label{b3}
If $Dw_2(\C B)=[\R B]$ then
\begin{displaymath}
\chi(\R B)\equiv\sigma(\C B)\pmod{8}
\end{displaymath}
\end{prop}

{}From \ref{b1} one can deduce a new congruence for complex orientations
of curves on a hyperboloid.
To formulate these congruences we shall use the integral calculates based
on Euler characteristic (see O.Ya.Viro\cite{V2}).
Let $H$ be a hyperboloid, $A$ be a nonsingular real curve of type I and
bidegree $(d,r)$ in $H$.
If curve $\R A$ supplied with the complex orientation realizes zero
in $H_1(\R H;\Z_4)$
then every fixed component of $\R H\setminus\R A$ determines function
$ind:\R H\setminus\R A\rightarrow\Z_4$ -- the linking number of point
and curve $\R A$ supplied with the complex orientation.
Let $\R A$ have components nonshrinking in $\R H$ and these components
realize the element with coordinates $(s,t),s,s\geq 0$ in standard
basis of $H_1(\R H)\approx H_1(\R P^1\times\R P^1)$.
Let $\R A$ supplied with the complex orientation realize class
$l'(s,t),l'\in\Z$
in $H_1(\R H)$.

\begin{prop}
\label{b4}
If $l'\equiv 0\pmod{4}, ds+rt\equiv 0\pmod{4}$ then
\begin{displaymath}
\int_{\R H}ind^2d\chi\equiv\frac{dr}{2}\pmod{8}
\end{displaymath}
\end{prop}

\begin{prop}
\label{b5}
If $l'\equiv 4\pmod{8}, ds+rt\equiv 2\pmod{4}$ then
\begin{displaymath}
\int_{\R H}ind^2d\chi\equiv\frac{dr}{2}+4\pmod{8}
\end{displaymath}
\end{prop}

\begin{prop}
\label{b6}
If $l'\equiv 0\pmod{8}$ then
\begin{displaymath}
\int_{\R H}ind^2d\chi\equiv\frac{dr}{2}\pmod{8}
\end{displaymath}
\end{prop}

\begin{prop}
\label{b7}
If $l'\equiv 0\pmod{8}, ds+rt\equiv 0\pmod{4}$ then
\begin{displaymath}
\int_{\R H}ind^2d\chi\equiv\frac{dr}{2}\pmod{16}
\end{displaymath}
\end{prop}

Let us now consider the relative case.
Let $\R A$ be a nonsingular real algebraic curve in $\R B$,
$\C A$ be the complexification of $\R A$.
Let $e_A$ denote the normal Euler number of $\C A$ in $\C B$.

\begin{prop}
\label{b8}
If $Dw_2(\C B)+[\R B]+[\C A]=0$ then there exists a natural separation of
$\R B\setminus\R A$ into two surfaces $B_1$ and $B_2$, $\partial B_1=\partial
B_2=\R A$.
This separation is determined by the condition that $\C A/conj\cup B_j$
is a characteristic surface in $\C B/conj$.
There is Guillou-Marin form $q_j$ on $H_1(\C A/conj\cup B_j;\Z_2)$
and there is congruence
\begin{displaymath}
\chi(B_j)\equiv\frac{e_A}{4}+\frac{\chi(\R B)-\sigma(\C
B)}{4}+\beta(q_j)\pmod{8}
\end{displaymath}
\end{prop}

\begin{prop}
\label{b9}
If $Dw_2(\C B)+[\R B]+[\C A]=0$ then difference $q_j(x)-q(x),x\in
H_1(B_j;\Z_2)$
is equal to the linking number of $x$ and $\C A$ in $\C B$, where $q$ is the
form from \ref{b1} and $q_j$ is the form from \ref{b8}
\end{prop}
The following result is an application of \ref{b8} and \ref{b9} to curves
on a hyperboloid.
Let $B_+$ be one of two parts of hyperboloid $\R H$ bounded by curve $\R A$
in the case when $d\equiv r\equiv 0\pmod{2}$.

\begin{prop}
\label{b10}
Let $d\equiv r\equiv 0\pmod{2}$ and $\frac{d}{2}t+\frac{r}{2}s+s+t\equiv
1\pmod{2}$.
\begin{itemize}
\begin{description}
\item[a)] If $A$ is an M-curve then $\chi(B_+)\equiv\frac{dr}{2}\pmod{8}$
\item[b)] If $A$ is an (M-1)-curve then $\chi(B_+)\equiv\frac{dr}{2}\pm
1\pmod{8}$
\item[c)] If $A$ is an (M-2)-curve and $\chi(B_+)\equiv\frac{dr}{2}+4\pmod{8}$
then $A$ is of type I
\item[d)] If $A$ is of type I then $\chi(B_+)\equiv 0\pmod{4}$
\end{description}
\end{itemize}
\end{prop}

Points a) and b) of \ref{b10} in the case when $\frac{d}{2}t+\frac{r}{2}s\equiv
0\pmod{2}, s+t\equiv 1\pmod{2}$ were proved by S.Matsuoka \cite{Ma1} in another
way (using 2-sheeted branched coverings of hyperboloid).
Point d) of \ref{b10} is a corollary of the modification of Rokhlin formula
of complex orientations for modulo 4 case.

The following theorem is a generalization of Rokhlin and
Kharlamov-Gudkov-Krakhnov congruences for curves on surfaces, because if $B$ is
an M-surface then
$Dw_2(\C B)=[\R B]$.
Let surface $\C B$ be the transversal intersection of hypersurfaces in $\C P^q$
defined by equations $P_j(x_0,\ldots,x_q)=0,j=1,\ldots,s-1$.
Let curve $\C A$ be the transversal intersection of $\C B$ and the
hypersurface in $\C P^q$ defined by equation $P_s(x_0,\ldots,x_q)=0$,
where $P_j$ is a real homogeneous polynomial, $deg P_j=m_j$, $q=s+1$.
Let $d$ denote $rk (in_*: H_1(B_+;\Z_2)\rightarrow H_1(\R B;\Z_2))$.

\begin{prop}
\label{b11}
Let $B$ be of type I $abs$ in the case when $\sum^{s-1}_{j=1}=s\pmod{2}$
or type I $rel$ in the case when $\sum^{s-1}_{j=1}=s+1\pmod{2}$
(for the definitions of type of surfaces see \cite{V}).
Let $B_+$ lie completely in one surface of complex separation.
Let homomorphism $inj_*:H_1(\R A;\Z_2)\rightarrow H_1(\R B;Z_2)$
be equal to zero.
Let $A$ be an (M-$k$)-curve.
Let $m_s$ be even and let, in the case when $m_s\equiv 2\pmod{4}$,
surface $B_-$ is contractible in $\R P^q$.
\begin{itemize}
\begin{description}
\item[a)] If $d+k=0$ then
$\chi(B_+)\equiv\frac{m_1\ldots m_{s-1}m_s^2}{4}\pmod{8}$
\item[b)] If $d+k=1$ then
$\chi(B_+)\equiv\frac{m_1\ldots m_{s-1}m_s^2}{4}\pm 1\pmod{8}$
\item[c)] If $d+k=2$ and
$\chi(B_+)\equiv\frac{m_1\ldots m_{s-1}m_s^2}{4}+4\pmod{8}$
then $A$ is of type I and $B_+$ is orientable
\item[d)] If $A$ is of type I and $B_+$ is orientable then
$\chi(B_+)\equiv\frac{m_1\ldots m_{s-1}m_s^2}{4}\pmod{4}$
\end{description}
\end{itemize}
\end{prop}

Next theorem allows to calculate the $\Z_2$-Netsvetaev class (see \cite{V})
for regular complete intersections.
It is easy to see that $\C B$ is spin manifold iff $\sum_j (m_j-1)\equiv 1
\pmod{2}$.
Denote by $\alpha\in H^1(\R B;\Z_2)$ the class induced under the inclusion map
$:\R B\rightarrow\R P^q$ from the only nontrivial element of $H^1(\R
P^q;\Z_2)$.

\begin{prop}
\label{b12}
If $\sum_j (m_j-1)\equiv 1\pmod{2}$ then the $\Z_2$-Netsvetaev class is
well-defined and equal to
\begin{displaymath}
\frac{1}{2} (1+\sum_j (m_j-1))\alpha \in H^1(\R B;\Z_2)
\end{displaymath}
\end{prop}

The author is indebted to O.Ya.Viro for the consultations.

\end{document}